\documentclass[11pt,leqno]{article}
\usepackage{amsmath, amssymb, amsfonts, amsthm,amscd}

\swapnumbers
\theoremstyle{plain}
\newtheorem{theorem}{Theorem}[section]
\newtheorem{lemma}[theorem]{Lemma}
\newtheorem{proposition}[theorem]{Proposition}
\newtheorem{corollary}[theorem]{Corollary}

\theoremstyle{definition}
\newtheorem{definition}[theorem]{Definition}
\newtheorem{remark}[theorem]{Remark}

\numberwithin{equation}{theorem}

 \newcommand{\Cont}{\operatorname{Cont}}
 \newcommand{\Imo}{\operatorname{Im}}  
  \newcommand{\om}{\omega}    \newcommand{\Om}{\Omega}
  \newcommand{\sig}{\sigma}   
  \newcommand{\al}{\alpha}
  \newcommand{\del}{\delta}   
  \newcommand{\gam}{\gamma}   
  \newcommand{\lam}{\lambda} 
  \newcommand{\veps}{\varepsilon}

  \def\b1{\text{\large 1}}  
  \def\wed{\wedge}
  \def\obulltimes{\overset{\bullet}{\otimes}}
  \def\half{\frac{\textup{\footnotesize 1}}{\textup{\footnotesize 2}}}
  \def\fourth{\frac{\textup{\footnotesize 1}}{\textup{\footnotesize 4}}}
  \def\simto{\overset{\sim}{\longrightarrow}}
  \def\tosim{\underset{\sim}{\longrightarrow}}
  \def\ip<#1>{\langle#1\rangle}   

  \newcommand{\ad}{\operatorname{ad}}
  \newcommand{\Ad}{\operatorname{Ad}}
  \newcommand{\hor}{\operatorname{hor}} 
  \newcommand{\Cl}{\operatorname{Cl}}
  \newcommand{\End}{\operatorname{End}}
    
  \newcommand{\Ind}{\operatorname{Ind}}  
  \newcommand{\Ker}{\operatorname{Ker}}
  
   \newcommand{\Exp}{\operatorname{Exp}}
 
\newcommand{\beqn}{\begin{equation}}
\newcommand{\eeqn}{\end{equation}}

 \newcommand{\fg}{\mathfrak{g}}
 \newcommand{\fp}{\mathfrak{p}}
 \newcommand{\fr}{\mathfrak{r}}

  \def\ghat{\hat{\mathfrak{g}}}

\newcommand{\bc}{\mathbb{C}}

\newcommand{\bz}{\mathbb{Z}}

 \newcommand{\cd}{\mathcal{D}}
 \newcommand{\ch}{\mathcal{H}}

 \newcommand{\cq}{\mathcal{Q}}

 \newcommand{\cw}{\mathcal{W}}
 
\newcommand{\ca}{\mathcal{A}}
 \newcommand{\cf}{\mathcal{F}}

\begin{document}

  \title{Induction Functor in Non-commutative 
Equivariant Cohomology and Dirac Cohomology}
  \author{Shrawan Kumar\\
   Department of Mathematics\\
   University of North Carolina\\
   Chapel Hill, N.C.  27599-3250}

  \maketitle

{\it Dedicated to Bertram Kostant on his seventy fifth birthday}

  \begin{abstract}
 The aim of this paper is to put some recent results of Huang-Pand\u{z}i\'c 
(conjectured by Vogan) and Kostant on Dirac cohomology in a broader 
perspective.  This is achieved by introducing an induction functor in the 
noncommutative equivariant cohomology. In this context, the results of 
Huang-Pand\u{z}i\'c and Kostant are interpreted as special cases (corresponding to the manifold being a point) of 
more general results on noncommutative equivariant cohomology introduced
 by Alekseev-Meinrenken. 
  \end{abstract}

  \section*{Introduction}

Let $G$ be a (not necessarily connected)  real Lie group and let $R$ be a 
 closed subgroup with their complexified Lie algebras 
$\fg$ and $\fr$ respectively.  We assume that there exists a $G$-invariant 
nondegenerate symmetric bilinear form $B_\fg$ on $\fg$ such that 
$B_{\fg_{|\fr}}$ 
is again nondegenerate.  We will impose this restriction on $G$ and $R$ 
throughout the paper.  Let $\fp$ be the orthocomplement $\fr^\perp$ of 
$\fr$ in $\fg$.  Then $B_{\fg_{|\fr}}$ being nondegenerate we have $\fg = 
\fr\oplus\fp$ and, moreover, $B_{\fg_{|\fp}}$ is again nondegenerate.  
Further $\fp$ is $R$-stable under the adjoint action.  For example, any 
compact  Lie group $G$ and a  closed subgroup $R$ 
satisfies the above restriction.

Let $M$ be a smooth $R$-manifold.  Then the deRham complex $\Om (M)$ of 
$M$ is canonically a $\bz_+$-graded (and hence $\bz/(2)$-graded) 
$R$-differential algebra.  We will only consider 
 $\bz/(2)$-graded spaces, algebras etc., so,  in the sequel, by graded we will mean  $\bz/(2)$-graded. We define a certain induction functor in 
noncommutative equivariant cohomology which associates to the 
$R$-differential algebra $\Om (M)$ a differential graded algebra $\Ind_{G/R}(\Om 
(M))$.  By definition,
  \[
\Ind_{G/R}(\Om (M)) = (\cw (\fg )\otimes\Om (M))_R, 
  \]
where $\cw (\fg ) := U(\fg )\otimes\Cl (\fg )$ is the noncommutative Weil 
algebra (cf. \S 1), $U(\fg )$ is the enveloping algebra, $\Cl (\fg )$ is 
the Clifford algebra of $\fg$ with respect to the form $B_\fg$, and the 
subscript $R$ refers to the subspace of `$R$-basic' elements (cf. \S1).  
The differential graded algebra structure on $\Ind_{G/R}(\Om 
(M))$ is the restriction of the tensor product differential graded algebra structure on $\cw (\fg )\otimes\Om (M)$.  We prove 
that the differential graded algebra $\Ind_{G/R}(\Om (M))$ is 
canonically isomorphic (as a differential graded algebra) with the 
non-commutative $G$-equivariant Cartan model $(U(\fg )\obulltimes \Om 
(M_G))^G$ of the $G$ manifold $M_G := G \times^R M$ (cf. Theorem 2.2).  
From this isomorphism,  we 
obtain (as an immediate corollary) that the cohomology $H(\Ind_{G/R}(\Om 
(M))$ is canonically isomorphic with the noncommutative $G$-equivariant 
cohomology $\ch_G(M_G)$ as graded algebras.

We use the above isomorphism to construct a functorial graded linear 
cochain map $\Phi_M: \Ind_{G/R}(\Om (M)) \to (U(\fr )\obulltimes \Om 
(M))^R$, where the latter is the noncommutative $R$-equivariant Cartan 
model of the $R$-manifold $M$.  Further, we show that $\Phi_M$ induces an 
algebra isomorphism in cohomology, even though, in general, $\Phi_M$ by 
itself is {\em not} an algebra homomorphism.  As a corollary, we obtain a 
functorial graded algebra isomorphism $\ch_G(M_G) \simeq \ch_R(M)$.

We now specialize the above results to the case when $M$ is the one point 
manifold $M^o$ to obtain some important recent results of 
Huang-Pand\u{z}i\'c and Kostant on Dirac cohomology ([Ko$_3$], [HP]).  In 
more detail, taking $M=M^o$, 
  \[
\Ind_{G/R}(\Om (M^o)) \cong (U(\fg )\otimes\Cl (\fp ))^R.
  \]
We show that the differential $d$ on $\Ind_{G/R}(\Om (M^o))$ corresponds 
under the above isomorphism with the differential $\ad \cd^\fp$ on the right 
side introduced by Kostant, where $\cd^\fp\in (U(\fg )\otimes\Cl (\fp ))^R$ is 
his remarkable cubic Dirac operator $\sum p_\ell \otimes q_\ell + 1\otimes \gamma_\fp$, 
where $\{ p_\ell\}_\ell$ is any basis of $\fp$ and $\{ q_\ell\}_\ell$ is 
the dual basis with respect to $B_{\fg_{|\fp}}$  and $\gamma_\fp$ is the Cartan 
element in $\wed^3(\fp )$ under the standard identification $\wed (\fp 
)\simeq \Cl (\fp )$.  Recall that in the case when $R$ is a maximal 
compact subgroup of reductive $G$, then $\gamma_\fp =0$ and the operator $\cd^\fp$ reduces to 
the Dirac operator considered by Vogan in defining his Dirac cohomology.  
Thus our theorem in the case $M=M^o$ gives that
  \beqn
H\bigl( (U(\fg )\otimes\Cl (\fp ))^R, \ad \cd^\fp\bigr) \simeq \ch_R(M^o) 
\simeq U(\fr )^R,  \tag{$*$}
  \eeqn
which was proved by Huang-Pand\u{z}i\'c in the case when $R$ is a maximal compact subgroup of a connected reductive $G$ and by Kostant in the general connected reductive case, i.e., when
$G$ and $R$ are connected and reductive (and of course $B_{\fg_{|\fr}}$ 
is  nondegenerate).  In fact from ($*$) one obtains the 
decomposition
  \[
\Ker (\ad\cd^\fp ) = \xi (Z(R )) \oplus \text{Image} (\ad \cd^\fp ), 
  \]
where the homomorphism $\xi : U(\fr )\to U(\fg )\otimes\Cl (\fp )$ is induced from the 
adjoint action of $\fr$ on $\fg$ and $\fp$ and $Z(R)$ is the subalgebra  of 
$R$-invariants $U(\fr )^R$.  Also the isomorphism ($\ast$) 
gives rise to an algebra homomorphism $\eta_R : Z(G ) \to Z(R )$. We show that, from the general properties of the Duflo isomorphism,  $\eta_R$ 
 is the unique homomorphism making the following diagram commutative:
  \[  \begin{CD}  
Z(G ) @>{\eta_R}>> Z(R ) \\
@V{H_\fg}VV @VV{H_\fg}V \\
S(\fg )^G @>>> S(\fr )^R , 
  \end{CD}  \]
where  $H_\fg$ is the Harish-Chandra 
isomorphism and the bottom horizontal map is induced by the orthogonal 
projection $\fg\to\fr$. 

\vskip2ex
\noindent
{\bf Acknowledgements.} It is my pleasure to thank S. Sahi who brought to my attention (and also explained) the work of Huang-Pand\u{z}i\'c [HP]. The work was partially supported from the NSF grant DMS 0070679.

  \section{Review of non-commutative 
equivariant cohomology after  Alekseev-Meinrenken}

Unless otherwise explicitly stated, by vector spaces we mean complex vector 
spaces and by linear maps complex linear maps.

  Let $G$ be a  (not necessarily connected) real Lie group with complexified
 Lie algebra $\fg$.  We assume that  $\fg$ has  a nondegenerate
symmetric $G$-invariant bilinear form $B_\fg$ on $\fg$, often denoted as
$\ip<\, ,\, >$. Define the
$\bz$-graded super-Lie algebra $\ghat^*$ as follows.  As a vector space,
  \[
\ghat^{-1} = \ghat^0 = \fg , \; \ghat^1 = \bc ,\; \ghat^n =0\text{ if } 
n\neq -1,0,1.
  \]
For $a\in\fg$, the corresponding element in $\ghat^{-1}$ (resp. $\ghat^0$) 
will be denoted by $i_a$ (resp. $L_a$) and they represent `contraction' 
and `Lie derivation' respectively.  We denote the generator of $\ghat^1$ by 
$d$.  The bracket relations in $\ghat^*$ are defined by (for 
$a,b\in\fg$): 
  \begin{align*}
[L_a, i_b] &= i_{[a,b]} \\
[L_a, L_b] &= L_{[a,b]} \\
[i_a, d] &= L_a \,.
  \end{align*}

By a {\em super-space} we mean a $\bz /(2)$-graded space. Any $\bz $-graded
space of course has a canonical $\bz /(2)$-grading by even and odd components.

 Recall that a  $G$-{\em differential space} is a super-space $B$ which is a 
Fr\'echet space, together with a graded smooth  action of $G$ on $B$ and a  
super-Lie algebra homomorphism $\theta : \ghat^* \to \End_{\Cont} B$,
where   $\End_{\Cont} B$ denotes the continuous linear endomorphisms of $B$. Moreover, we assume that the action of $G$ commutes with $d$ , $L_a$ is the derivative of the $G$-action and 
$$ gi_ag^{-1}=i_{g\cdot a} \,\,\text{for all }\, g\in G \,, a\in \fg.$$
The {\em 
horizontal subspace} $B_{\hor}$ is the space annihilated by $\ghat^{-1}$, 
the {\em invariant subspace} $B^G$ is the subspace invariant under $G
$ 
and the space $B_G$ of {\em basic elements} is the intersection $B_{\hor} 
\cap B^G$.

A {\em $G$-differential algebra} is a super-algebra $B$ together with the 
structure of a $G$-differential space on $B$ such that $\theta$ takes 
values in the super-derivations of $B$ and the action of $G$ on $B$ is via algebra automorphisms.

For a smooth $G$-manifold $M$,  $\Omega (M)$ with the Fr\'echet topology
provides the most important class 
of examples of $G$-differential algebras.

A {\em homomorphism} between $G$-differential spaces (resp. algebras) 
$(B_1,\theta_1)$ and $(B_2, \theta_2)$ is a continuous homomorphism of super-spaces 
(resp. algebras) $\phi: B_1\to B_2$ such that
  \[\phi(g\cdot b)=g\cdot \phi(b),\,\,
\phi (\theta_1(x)b) = \theta_2(x)\,\phi (b), \;\text{ for $g\in G, x\in\ghat^*$ and 
$b\in B_1$}. 
  \]

There is the (classical) {\em Weil algebra} $W (\fg ) := S(\fg^*)\otimes\wed 
(\fg^*)$ with the tensor product algebra structure.  This is a 
$G$-differential algebra under the $\bz_+$-grading
  \[
W(\fg )^n = \bigoplus_{k\geq 0} S^k(\fg^*)\otimes\wed^{n-2k}(\fg^*).
  \]
The action of $G$ is via the coadjoint action. The operators $L_a$ come from the coadjoint action of $\fg$ on $S(\fg^*)$ 
and $\wed (\fg^*)$.  The contraction operator $i_a$ on $W(\fg )$ is 
defined as $I_{S(\fg^*)}\otimes i'_a$, $i'_a$ being the standard 
contraction operator on $\wed (\fg^*)$.  The differential $d$ on $W(\fg )$ 
is the unique super-derivation satisfying (for any $f\in\fg^*$)
  \[
d(\b1\otimes e_f) = \b1\otimes d_\wed e_f + s_f\otimes\b1 ,
  \]
where $e_f$ (resp. $s_f$) is the element $f$ considered as an element of
$\fg^*\subset\wed (\fg^*)$ (resp. $\fg^*\subset S(\fg^*)$) and $d_\wed : 
\wed (\fg^*) \to \wed (\fg^*)$ is the standard Koszul 
differential.

We now recall the definition of the {\em non-commutative Weil algebra} 
$\cw (\fg )$ from [AM], which is a $G$-differential algebra.  Recall that 
the Clifford algebra $\Cl (\fg )$ of $\fg$ with respect to the 
bilinear form $B_\fg$ is the quotient of the tensor 
algebra $T(\fg )$ of $\fg$ by the two-sided ideal generated by $2\, 
x\otimes x - \ip<x,x>, \; x\in\fg$.  As a super-space it is defined as 
the tensor product of algebras
  \[
\cw (\fg ) := U(\fg )\otimes \Cl (\fg ), 
  \]
where the $\bz /(2)$-grading on $\cw (\fg )$ comes from the standard 
grading on the Clifford algebra $\Cl (\fg )$, and $U(\fg )$ is the enveloping algebra of $\fg$ placed in the even degree part.
Both of $U(\fg )$ and $\Cl (\fg )$ are $G$-modules under the adjoint 
action and so is their tensor product.  For $a\in\fg$, let $L_a$ be the 
adjoint action on $\cw (\fg )$.

Recall that there is a vector space isomorphism given by the symbol map
  \[
\sig : \Cl (\fg ) \to \wed (\fg ),
  \]
where $\sig^{-1}$ is induced from the standard projection map $T(\fg) \to 
 \Cl (\fg )$ under the identification of $\wed (\fg )$ with the skew-symmetric
tensors in $T(\fg)$. {\it From now on we will  identify $\Cl (\fg )$ with $\wed (\fg )$ (via the symbol 
map) as a vector space.}  Under this identification, we will denote the 
product in $\wed (\fg )$ by $\odot$, i.e.,
  \[
x\odot y = \sig \bigl(\sig^{-1}(x)\cdot \sig^{-1}(y)\bigl) , \;\text{ for 
} x,y\in\wed (\fg ).
  \]
The exterior product in $\wed (\fg )$ will be denoted by $x\wed y$.   
Recall that $\wed (\fg )$ admits the contraction operator $\bar{i}_a$ (for 
$a\in\fg$) which is a super-derivation induced from the operator
  \[
\bar{i}_a b = \ip<a,b>, \;\text{ for } b\in\fg .
  \]
Then, by [AM, Lemma 3.1], explicitly the product $\odot$ in $\wed (\fg )$ 
is given by (for $\om, \eta \in\wed (\fg )$):
  \[
\om\odot\eta = \mu \Bigl(\Exp \Bigl( -\half \sum_k\, i^1_{a_k} 
i^2_{b_k}\Bigr) (\om\otimes\eta )\Bigr) ,
  \]
where $\{ a_k\}_k$ is any basis of $\fg$ and $\{ b_k\}_k$ is the dual 
basis $\ip<a_k,b_\ell> = \del_{k,\ell}$, $\mu : \wed (\fg )\otimes \wed 
(\fg ) \to \wed (\fg )$ is the standard wedge product, 
$i^1_{a_k}(\om\otimes\eta ) := (\bar{i}_{a_k}\om )\otimes\eta$, 
$i^2_{a_k}(\om\otimes\eta ) = \om\otimes \bar{i}_{a_k}\eta$.

 Define the 
operator $i_a$, $a\in \fg$,  on $\cw (\fg )$ by
  \[
i_a = I_{U(\fg )}\otimes \bar{i}_a .
  \]
For $a\in\fg$, let $u_a$ (resp. $c_a$) be the corresponding element in 
$\fg\subset U(\fg )$ (resp. $\fg\subset\Cl (\fg )$).
We also think of $u_a$ (resp. $c_a$) as the element $u_a\otimes\b1$  
(resp. $\b1 \otimes c_a$) of $\cw (\fg )$.

Finally, we define the differential $d: \cw (\fg ) \to \cw (\fg )$ as the 
commutator
   \[  dx = \ad \cd (x),     \]
where $\cd\in\cw (\fg )$ is defined by
  \[  \cd = \sum_k u_{a_k} c_{b_k} - \b1\otimes \gam ,  \] $\gam = \gam_\fg \in \wed^3(\fg )^G$ is the $G$-invariant 
element   (so called the {\it Cartan element}) defined by
  \[
\gam (a,b,c) = \ip<a, [b,c]>, \;\text{ for } a,b,c\in\fg 
  \]
under the identification 
$\wed (\fg ) \simeq \wed (\fg^*)$ induced from the form $\ip<\, ,\,>$, and 
$\ad \cd$ is the super-adjoint action defined by $\ad \cd (x)=\cd x-x\cd$,
for $x\in \cw (\fg )^{\text even}$, and  $\ad \cd (x)=\cd x+x\cd$,
for $x\in \cw (\fg )^{\text odd}$. 

Then $\cw (\fg )$ with the above operators $i_a$, $L_a$, $d=\ad \cd$ and the 
adjoint action of $G$  
becomes a $G$-differential algebra called the {\it non-commutative Weil 
algebra}.

By [AM, Proposition 3.7 and the equation (3)], $d$ is given by the formula
  \begin{align*}
d(x\otimes\om) = & -\ad u_{a_k} ( x) \otimes c_{b_k}\wed\om - 
\Biggl(\frac{u_{a_k}x+xu_{a_k}}{2}\Biggr) \otimes \bar{i}_{b_k}\om \\
&\; + x\otimes d_\wed \om +\frac{1}{4} x\otimes \bar{i}_\gam \om ,\;\text{for } x\in U(\fg 
), 
\om\in\wed (\fg ),
  \end{align*}
where $d_\wed$ is the Koszul differential on $\wed (\fg )$ of degree +1 
under the identification $\wed (\fg ) \simeq \wed (\fg^*)$.

Recall [C] that the {\em $G$-equivariant cohomology} $H_G(B)$ of a $G$-differential 
algebra $B$ is by definition the cohomology of the basic subalgebra 
$(W(\fg )\otimes B)_G$ of the tensor product $G$-differential algebra 
$W(\fg)\otimes B$ under the tensor product differential $d(x\otimes y) = 
dx\otimes y + (-1)^{\deg x}\, x\otimes dy$, for $x\in W(\fg )$ and $y\in 
B$.

Similarly, following [AM], the {\it non-commutative $G$-equivariant cohomolgy} 
$\ch_G(B)$ is the cohomology of the basic subalgebra $(\cw (\fg )\otimes 
B)_G$ under the tensor product differential. Then, clearly, $\ch_G(B)$ is
a super-algebra. 

  \begin{proposition}  For any $G$-differential algebra $B$, the projection map 
$\theta : W(\fg )\otimes B\to S(\fg^*)\otimes B$, induced from the 
standard augmentation map $\wed (\fg^*)\to\bc$, induces an algebra 
isomorphism (again denoted by)
  \[
\theta : \bigl( W(\fg )\otimes B\bigr)_G \simto \bigl( S(\fg^*)\otimes 
B\bigr)^G .
  \]

Under the above isomorphism, the differential $d$ corresponds to the 
differential $d_G$ on $\bigl( S(\fg^*)\otimes B\bigr)^G$ given as follows:
  \[
d_G = I_{S(\fg^*)} \otimes d - \sum_k s_{a^*_k}\otimes i_{a_k} ,
  \]
where $\{ a_k\}_k$ is a basis of $\fg$ and $\{ a^*_k\}_k$ is the dual 
basis of $\fg^*$ and $s_{a^*_k}$ denotes the operator acting on $S(\fg^*)$
via the multiplication by $a^*_k$. 
  \end{proposition}

Similarly, we have the following proposition from [AM, \S4.2].

  \begin{proposition}  For any $G$-differential algebra $B$, the projection map
  \[
\Theta : \cw (\fg )\otimes B \to U(\fg )\otimes B,
  \]
induced from the standard augmentation map $\wed (\fg 
)\to\bc$, induces a vector space (but not in general algebra) isomorphism
  \[
\Theta : (\cw (\fg )\otimes B)_G \simto (U(\fg )\otimes B)^G .
  \]
To distinguish, let $(U(\fg )\obulltimes B)^G$ denote the vector space 
$(U(\fg )\otimes B)^G$ with the new product $\odot$ making $\Theta$ an algebra 
isomorphism.

Under the above isomorphism, the differential $d$ corresponds to the 
differential
  \beqn
d_G = I_{U(\fg )}\otimes d - \half \sum_k \bigl( u^L_{a_k} + 
u^R_{a_k}\bigr) \otimes i_{b_k} + \fourth\, I_{U(\fg )}\otimes 
i_\gam ,
  \eeqn
where $\gam\in\wed^3(\fg )^G$ is defined earlier, $u^L_{a_k}$ (resp. 
$u^R_{a_k}$) denotes the left (resp. right) multiplication in $U(\fg )$ by 
$u_{a_k}$ and $\{ a_k\}_k, \{ b_k\}_k$ are dual bases of $\fg$.

 \end{proposition}

By [AM, Proposition 4.3], explicitly the multiplication $\odot$ in $(U(\fg 
)\otimes B)^G$ is given as the restriction of the multiplication (again 
denoted by) $\odot$ in $U(\fg )\otimes B$ defined as follows.  For $x,y\in 
U(\fg )$, $b_1, b_2\in B$, 
  \beqn
(x\otimes b_1)\odot (y\otimes b_2) = xy\otimes\mu \Bigl( \Exp\Bigl( -\half 
\sum_k i^1_{a_k} i^2_{b_k}\Bigr) (b_1\otimes b_2)\Bigr) ,
   \eeqn
where $i^1_{a_k}$ and  
$i^2_{b_k}$  are the contraction operators on $B\otimes B$ with respect
to the first and second factors respectively and  
$\mu : B\otimes B \to B$ is the multiplication map.

As in [AM], there exists a $G$-module  isomorphism (depending only on $\fg$), called the 
{\it quantization map},
  \[  \cq =\cq_\fg : W(\fg ) \to \cw (\fg )   \]
which intertwines all the operators $L_a$, $i_a$ and $d$.  
$\cq_{|S(\fg^*)}$ is the composite of the isomorphisms
  \[
S(\fg^*) \to S(\fg ) \overset{D_\fg}{\longrightarrow} U(\fg ),
  \]
where the first map is the algebra isomorphism induced from the 
isomorphism $\fg^* \to \fg$ (coming from $\ip<\, ,\,>$) and $D_\fg$ is the 
Duflo isomorphism [D]. (Recall that  $D_\fg$ is only a linear isomorphism from 
 $S(\fg )$ to $ U(\fg )$ but is an algebra isomorphism  restricted to $S(\fg 
)^\fg$ onto the center $U(\fg )^\fg$. Moreover, it maps isomorphically
 $S(\fg 
)^G$ onto  $U(\fg )^G$.) Also recall that, for $a\in S^p(\fg)$,
$D_\fg (a)=\Sigma(a)$ (mod $U(\fg)^{p-1}$), where $U(\fg)^{p}$
is the standard filtration of $U(\fg)$ and $\Sigma: S(\fg) \to U(\fg)$ is the standard symmetrization map. Also $\cq_{\vert\wed (\fg^*)}$ is the 
isomorphism (induced from $\ip<\, ,\,>$)
  \[   \wed (\fg^*) \simto \wed (\fg ).   \]
Of course, as earlier, we have identified $\Cl (\fg )$ with $\wed (\fg )$ 
via the symbol map $\sig$.

However, $\cq \neq \cq_{|S(\fg^*)} \otimes \cq_{|\wed (\fg^*)}$ in general.

  \begin{theorem} {[AM, Theorem 7.1]}  For any $G$-differential algebra 
$B$, 
the cochain map $\cq\otimes I_B: W(\fg )\otimes B \to \cw (\fg )\otimes B$ 
induces an algebra isomorphism in cohomology
  \[
\cq^B : H_G(B) \simto \ch_G(B).
  \]
  \end{theorem}

 \begin{definition}
Let $\hat{\cq}^B_G = \hat{\cq}^B : (S(\fg^*)\otimes B)^G \to (U(\fg )\obulltimes 
B)^G$ be the unique map making the following diagram commutative:
  \[ \begin{CD}
(W(\fg )\otimes B)_G @>{\cq\otimes I_B}>{\sim}> (\cw(\fg )\otimes B)_G \\
 @VV{\theta}V @VV{\Theta}V \\
(S(\fg^*)\otimes B)^G @>{\hat{\cq}^B}>{\sim}> (U(\fg )\obulltimes B)^G.
  \end{CD}  \]
Then clearly $\hat{\cq}^B$ is a cochain isomorphism.  In general, 
$\hat{\cq}^B$ is not an algebra homomorphism.
  \end{definition}

  \section{An induction functor in non-commutative\\ equivariant 
cohomology}

Let $G$ be a  real (not necessarily connected) Lie group with complexified Lie 
algebra $\fg$ and let $R$ be a closed subgroup of $G$ with 
complexified Lie algebra $\fr$.
Assume that  $\fg$ admits a 
$G$-invariant nondegenerate symmetric bilinear form $B_\fg = \ip<\, 
,\,>$ such that ${B_\fg}_{|\fr}$ is 
nondegenerate.  We call such a pair of $(G,R)$ a {\it quadratic pair}. Thus we have the decomposition
  \[
\fg = \fr\oplus\fp , \quad \fp := \fr^\perp .
  \]
By the $G$-invariance of $B_\fg$, $\fp$ is $R$-stable under the 
adjoint action. Moreover, $ {B_\fg}_{|\fp}$ also is 
nondegenerate.

The following definition is influenced by the corresponding definition in the (commutative) equivariant cohomology given in [KV, Definition 32].

  \begin{definition}[Induction functor]  For an $R$-differential complex 
$B$, define the cochain complex
  \[
\Ind_{G/R}(B) = (\cw (\fg )\otimes B)_R
  \]
equipped with the standard tensor product differential
  \[
d(x\otimes y) = d_\cw x \otimes y + (-1)^{\deg x} x\otimes d_By,
  \quad x\in\cw (\fg ), y\in B,
  \]
where $d_\cw$ is the differential in $\cw (\fg )$ and $d_B$ is the 
differential in $B$.

Since $\cw (\fg )$ is a $G$ (in particular $R$) differential complex and 
$B$ is a $R$-differential complex, it is easy to see that indeed $d$ keeps 
the $R$-basic subspace of $\cw (\fg )\otimes B$ stable.

If $B$ is a $R$-differential   algebra, then $
\Ind_{G/R}(B)$ is a differential  algebra under the tensor product
super-algebra structure on  $\cw (\fg )\otimes B$. 

Let $M$ be a smooth real $R$-manifold and let $\Om (M)$ be the complexified 
deRham complex of $M$.    Consider the fiber product 
$G$-manifold $M_G := G\times^RM$, where $G$ acts on $M_G$ via the left 
multiplication on the first factor.
  \end{definition}

  \begin{theorem}  There exists a graded  algebra isomorphism
  \[
\psi_M : \bigl( U(\fg )\obulltimes \Om (M_G)\bigr)^G \simto \Ind_{G/R} 
(\Om (M)) 
  \]
commuting with the differentials, where 
$(U(\fg )\obulltimes \Om  
(M_G))^G$ 
is equipped with the Cartan differential $d_G$ (cf. Proposition 1.2).

Moreover, $\psi_M$ is functorial in the sense that for any $R$-equivariant 
smooth map $f: M\to N$, the following diagram is commutative:
  \[ \begin{CD}
(U(\fg )\obulltimes \Om (N_G))^G @>\psi_N>> \Ind_{G/R}(\Om (N)) \\
@V{I_{U(\fg )}\otimes f^*_G}VV @VV{I_{\cw (\fg )}\otimes f^*}V\\
(U(\fg )\obulltimes\Om (M_G))^G @>\psi_M>> \Ind_{G/R}(\Om (M)), 
  \end{CD}   \]
where $f^*: \Om (N) \to \Om (M)$ and  $f^*_G: \Om (N_G) \to \Om (M_G)$ 
are the induced maps from $f$.
  \end{theorem}

  \begin{proof}  Under the projection map $p: G\times M \to M_G$, 
we can identify
  \[
\Om (M_G) \subset \Om (G\times M).
  \]
For $\om\in\Om (M_G)$, by $\om (\b1 )$ we mean the evaluation of 
$\om$ at $\b1\times M$.  Thus $\om (\b1 )\in \wed (\fg^*)\otimes \Om (M)$.  
Under the identification $\wed (\fg^*)\simeq \wed (\fg )$ (induced from 
the bilinear form $B_\fg$), we can (and will) think of $\om (\b1 )\in \wed 
(\fg )\otimes\Om (M)$.  Thus we get the map
  \[
U(\fg )\otimes \Om (M_G) \to \cw (\fg )\otimes \Om (M), 
  \]
$x\otimes \om \mapsto x\otimes \om (\b1 )$.  Let $\psi_M$ be its 
restriction to $\bigl(U\bigl(\fg )\otimes\Om (M_G)\bigr)^G$. We need to show 
that  $ \psi_M \bigl(\bigl(U(\fg )\otimes\Om (M_G)\bigr)^G\bigr)\subset 
\Ind_{G/R}(\Om (M))$.

It is easy to see that (cf. [KV, Page 147])
  \[
\Om (M_G) = C^{\infty} \bigl( G, (\wed (\fg^*)\otimes \Om (M))_{\hor 
R}\bigr)^R ,
  \]
where $C^\infty (G,V)$ denotes the space of $C^\infty$-functions on $G$ 
with values in $V$; the $R$-invariants are taken with respect to the 
action of $R$ on $G$ via right multiplication, the given action of $R$ on 
$M$ and adjoint action on $\wed (\fg^*)$; the contraction $i_a$ 
($a\in\fr$) acting on $\wed (\fg^*)\otimes \Om (M)$ is the standard tensor 
product contraction.  Thus 
  \beqn
\bigl( U(\fg )\otimes\Om (M_G)\bigr)^G = \Bigl( C^\infty \bigl( G,(U(\fg 
)\otimes\wed (\fg^*)\otimes\Om (M))_{\hor R}\bigr)^R\Bigr)^G ,
  \eeqn
where $R$ acts trivially on $U(\fg )$; and the $G$-invariants are taken 
with respect to the left action of $G$ on $G$, the adjoint action on 
$U(\fg )$ and the trivial action on $\wed (\fg^*)\otimes\Om (M)$.

Take $\tilde{\al}\in \bigl( C^\infty (G, (U(\fg )\otimes \wed 
(\fg^*)\otimes\Om (M))_{\hor R})^R\bigr)^G$.  Then,
  \[
\tilde{\al} (gk^{-1}) = k\, \tilde{\al}(g), \;\text{ for } g\in G, k\in R.
  \]
Writing $\tilde{\al}(\b1 )=\sum x_i\otimes\om_i$, $x_i\in U(\fg )$, $\om_i\in \wed (\fg^*)\otimes\Om (M)$,  since 
$\tilde{\al}$ is $G$-invariant,
  \begin{align*}
\sum_i (\Ad (gk^{-1})\, x_i)\otimes \om_i &= k\,\cdot \sum_i (\Ad g\, 
x_i)\otimes \om_i \\
 &= \sum_i\, (\Ad g\, x_i)\otimes k\, \om_i .
  \end{align*}
Taken $g=k$ in the above identity, we get
  \[
\sum_i\, x_i\otimes\om_i = \sum_i\, (\Ad k\, x_i)\otimes k\om_i .
  \]
Thus $\psi_M(\tilde{\al}) \in \Ind_{G/R}(\Om (M))$.

We next show that $\psi_M$ is surjective onto $\Ind_{G/R}(\Om (M))$.  Take 
$\al = \sum_i x_i\otimes\om_i\in$ $\bigl( U(\fg )\otimes (\wed 
(\fg^*)\otimes\Om (M))\bigr)_R$ and define $\tilde{\al}\in (U(\fg 
)\otimes\Om (M_G))^G$, under the identification (1), by
  \[   
\tilde{\al}(g) = \sum_i \bigl(\Ad g\, x_i\bigr) \otimes \om_i, \;\; g\in G.
  \]
Clearly $\tilde{\al}$ is  $G$-invariant.  Further,
  \[
\tilde{\al} \in C^\infty \bigl( G, (U(\fg )\otimes \wed (\fg^*) \otimes \Om 
(M))_{\hor R}\bigr)^R.
    \]
To show this, it suffices to prove that for all $g\in G$ and $k\in R$, 
  \beqn
\sum_i\, (\Ad (gk^{-1}) x_i) \otimes \om_i = \sum\, \Ad g\, x_i \otimes 
k\om_i .
  \eeqn
But $\al$ being $R$-invariant, 
  \beqn
\sum_i\, x_i\otimes \om_i = \sum\, (\Ad k\,x_i)\otimes k\om_i, \;\;\text{for all 
}k\in R.
  \eeqn
Applying $gk^{-1}$ to (3) we get (2).

The injectivity of $\psi_M$ is clear from the $G$-invariance of any 
element in the domain of $\psi_M$.  Thus $\psi_M$ is a linear isomorphism.  
We next show that $\psi_M$ is a cochain map.

View $\wed (\fg^*)$ as the space of left invariant forms on $G$.  For 
$x\in U(\fg )$, $f\in C^\infty (G)$, $\om_1\in \wed (\fg^*)$ and 
$\om_2\in\Om (M)$, by (1.2.1) ($\bar{d}$ being deRham differentials on
$G$ and also on $M$), 
  \begin{align*}
d_G (x\otimes f\om_1\otimes\om_2) &= x 
\otimes\bar{d}f\wed\om_1\otimes\om_2 + x\otimes 
f d_\wed \om_1\otimes\om_2\\
&+ (-1)^{\text{deg} \om_1} x\otimes f\om_1\otimes\bar{d}\om_2 
- \half \sum_k\bigl( u_{a_k}x + xu_{a_k}\bigr)
  \otimes f\bigl( \bar{i}_{b_k}\om_1\bigr) \otimes \om_2
\\
&\quad  + 
\fourth x\otimes f(\bar{i}_\gam \om_1)\otimes\om_2 .
  \end{align*}
(In fact, to be precise, in the above we should have taken $\sum_j\,
x^j\otimes f^j\om_1^j\otimes\om_2^j \in \bigl( U(\fg )\otimes 
\Om (M_G)\bigr)^G$ instead of just a single term $x\otimes f\om_1\otimes\om_2$. But, for notational convenience, we take a single term.)

Moreover, for any $x\otimes f\om_1\otimes\om_2 \in$ $\bigl( U(\fg )\otimes 
\Om (M_G)\bigr)^G$, we get (for any $a_k\in\fg$)
  \beqn
( u_{a_k}x-x u_{a_k})\otimes f\om_1\otimes 
\om_2 = -x\otimes a_k(f)\, \om_1\otimes\om_2 .
  \eeqn
Thus,
  \begin{align}
\psi_Md_G(x\otimes f\om_1\otimes \om_2) &=  x\otimes (\bar{d}f) 
(1)\wed\om_1\otimes\om_2 + x\otimes f(1)\, d_\wed \om_1\otimes\om_2 
\notag\\
&\qquad + (-1)^{\deg \om_1} x\otimes f(1)\om_1\otimes \bar{d}\om_2 
\notag\\
&\qquad -\half \sum_k\bigl( u_{a_k}x+xu_{a_k}\bigr)\otimes 
f(1)\bigl( \bar{i}_{b_k}\om_1)\otimes\om_2  \notag\\
&\qquad +\fourth x\otimes f(1) (\bar{i}_\gam\om_1) \otimes\om_2 \notag\\
&=  -\sum_k \bigl( u_{a_k}x - xu_{a_k}\bigr) \otimes f(1) 
a^*_k\wed \om_1\otimes \om_2 \\
&\qquad + x\otimes f(1) \, d_\wed \om_1\otimes \om_2 
  + (-1)^{\deg \om_1} x \otimes f(1)\, \om_1\otimes\bar{d}\om_2  \notag\\
&\qquad -\half \sum_k\bigl( u_{a_k}x+xu_{a_k}\bigr)\otimes 
f(1)\bigl( \bar{i}_{b_k}\om_1)\otimes\om_2  \notag\\
&\qquad +\fourth x\otimes f(1) (\bar{i}_\gam\om_1) \otimes\om_2, \;\text{ 
by (4)}, \notag
   \end{align}
where $\{ a^*_k\}$ is the basis of $\fg^*$ dual to the basis $\{ a_k\}$ of 
$\fg$.

On the other hand, by the expression of $d_\cw$ given in Section 1, 
  \begin{align}
d(x\otimes &f(1)\, \om_1\otimes\om_2)\notag\\ 
&= f(1) d_\cw (x\otimes \om_1) \otimes \om_2 
+ (-1)^{\deg \om_1}f(1)\, x\otimes\om_1\otimes \bar{d}\om_2 \notag \\
&= f(1) \biggl( -\ad u_{a_k} (x)\otimes c_{b_k}\wed\om_1 - 
\bigl(\frac{u_{a_k}x+xu_{a_k}}{2}\bigr)\otimes \bar{i}_{b_k}\om_1 + 
x\otimes d_\wed \om_1\\ 
&\qquad\qquad\; + \fourth x\otimes \bar{i}_{\gam} \om_1\biggr) \otimes\om_2 
+ (-1)^{\deg \om_1} f(1)\, x\otimes \om_1\otimes \bar{d}\om_2 .\notag
  \end{align}
Comparing (5) and (6) we get that $\psi_M$ commutes with the 
differentials.

Finally, we show that $\psi_M$ is an algebra homomorphism.  Take two 
elements $u=x\otimes\sum_i f_i\om'_i\otimes\om''_i$ and $v = 
y\otimes\sum_j g_j\, \eta'_j\otimes\eta''_j$ in $U(\fg )\obulltimes \Om (M_G)$, 
where $f_i, g_j\in C^\infty (G\times M)$, $x,y\in U(\fg )$, $\om'_i, 
\eta'_j \in\wed (\fg^*)$ and $\om''_i, \eta''_j \in\Om (M)$.  Then,
by (1.2.2), 
  \begin{align*}
\psi_M(u\odot v) &= \psi_M\sum_{i,j} 
(-1)^{\text{deg} \eta'_j\text{deg} \om''_i}\biggl( xy\otimes f_ig_j\, \mu 
\biggl(\Exp \biggl( -\half \sum_k i^1_{a_k} i^2_{b_k}\biggr) 
(\om'_i\otimes\eta'_j)\biggr)\biggr) \otimes \om''_i\eta''_j \\
 &= \sum_{i,j} \,
(-1)^{\text{deg} \eta'_j\text{deg} \om''_i} xy\otimes\mu\biggl(\Exp\biggl( -\half \sum_k\, i^1_{a_k} 
i^2_{b_k}\biggr) (\om'_i\otimes\eta'_j)\biggr)\\
&\qquad\qquad \otimes f_i(1, -)g_j(1,-)\om''_i\eta''_j \\
 &= \sum_{i,j} \,
(-1)^{\text{deg} \eta'_j\text{deg} \om''_i} xy\otimes (\om'_i\odot \eta'_j) \otimes f_i(1,-)\, 
g_j(1,-)\om''_i\eta''_j \\
 &= \psi_M(u)\cdot \psi_M(v).
  \end{align*}
This completes the proof of the theorem.  
  \end{proof}

  \begin{corollary}  For a real $R$-manifold $M$, the cochain map $\psi_M$ 
induces a functorial graded  algebra isomorphism:
  \[
\psi^*_M : \ch_G(M_G) \simto H\bigr(\Ind_{G/R}(\Om (M))\bigr).
  \]
  \end{corollary}

   \begin{definition}  Let $M$ be a real $R$-manifold.  Define a 
cochain map $\Phi_M: \Ind_{G/R}(\Om (M))\to (U(\fr )\obulltimes \Om (M))^R,$
making the following diagram commutative:

  \[  \begin{CD} 
 (S (\fg^* ) \otimes\Om (M_G))^G @.  @>{\al_M}>> (S(\fr^*)\otimes\Om (M))^R 
\\
 \hat{\cq}^{M_G}_{G}{\swarrow}@. @.  @VV{\hat{\cq}^M_R}V \\
(U (\fg)\obulltimes\Om (M_G))^G @>{\psi_M}>{\sim}> 
\Ind_{G/R}(\Om (M)) @>{\Phi_M}>> 
(U(\fr )\obulltimes\Om (M))^R ,
   \end{CD}  \]  
where $\al_M$ is induced from the map
  \[
P\otimes\eta \mapsto P_{|\fr}\otimes \eta_{|1\times 
M},    
  \]
for $P\in S(\fg^*)$ and $\eta\in\Om (M_G)$, where 
$ \hat{\cq}^{M_G}_{G}$ (and $\hat\cq^M_R$)  is the map defined in 1.4.

In fact, define
  \[
\Phi_M := \hat{\cq}^M_R \circ \al_M \circ (\hat{\cq}^{M_G}_{G})^{-1} \circ 
(\psi_M)^{-1} .
  \]

Then, clearly it is a cochain map making the above diagram commutative.
Observe that, in general,  $\Phi_M$ is not an algebra homomorphism.

Further, $\Phi_M$ is functorial in the sense that for any $R$-equivariant 
smooth map $f: M\to N$, the following diagram is commutative:
  \[  \begin{CD}
\Ind_{G/R}(\Om (N)) @>{\Phi_N}>> (U(\fr )\obulltimes\Om (N))^R \\
 @VVV @VVV \\
\Ind_{G/R}(\Om (M)) @>{\Phi_M}>> (U(\fr )\obulltimes\Om (M))^R\,,
  \end{CD}  \]
where the vertical maps are induced canonically from the $R$-differential 
algebra homomorphism $f^*: \Om (N) \to \Om (M)$.
  \end{definition}

  \begin{theorem}    
For any $R$-manifold $M$, the cochain map $\Phi_M : \Ind_{G/R}(\Om (M))$ 
$\to (U(\fr )\obulltimes\Om (M))^R$ induces an algebra isomorphism in cohomology:
  \[
[\Phi_M] : H\bigl( \Ind_{G/R}(\Om (M))\bigr) \simto \ch_R(M).
  \]

Thus, by Corollary 2.3, we have a functorial algebra isomorphism
  \[
\ch_G(M_G) \simto \ch_R(M).
  \]
In particular, $\ch_G(G/R)\simto U(\fr)^R$.
  \end{theorem}

  \begin{proof}  In the first commutative diagram of (2.4), all the maps are 
cochain maps.  Moreover, all the cochain maps $ \hat{\cq}^{M_G}_{G}$, 
$\psi_M$, $\hat{\cq}^M_R$ are cochain isomorphisms.  So it 
suffices to prove that $\al_M$ induces an isomorphism in cohomology.  But 
this follows from [DV, Th\'eor\`eme 24].
  \end{proof}

Let $M$ be a $R$-manifold.
Consider the $R$-module isomorphism 
  \[
\cq_\fg\otimes I_{\Om (M)}:  W(\fg )\otimes\Om (M) \to \cw (\fg )\otimes \Om 
(M).
  \]
Since $\cq_\fg$ commutes with the operators $L_a, i_a$ ($a\in \fg$) and $d$; 
in particular, $\cq_\fg\otimes I_{\Om (M)}$ induces the map  
  \[
\cq_{G/R}^M : (W(\fg )\otimes\Om (M))_R \to \Ind_{G/R}(\Om (M))
  \]
commuting with differentials, where $W(\fg )\otimes\Om (M)$ is equipped 
with the standard tensor product $R$-differential algebra structure.

  \begin{lemma}  For any $R$-manifold $M$, the following diagram is 
commutative:
  \[ \begin{CD}
(W(\fg )\otimes\Om (M))_R @>{\hat{\al}_M}>> \bigl( S(\fr^*)\otimes\Om 
(M)\bigr)^R\\
@V{\cq^M_{G/R}}VV  @VV{\hat{\cq}^M_R}V \\
\Ind_{G/R}(\Om (M)) @>>{\Phi_M}> \bigl( U(\fr )\overset{\bullet}{\otimes} 
\Om (M)\bigr)^R ,
  \end{CD}  \]
where $\hat{\al}_M (P\otimes\om\otimes\eta ) = \veps (\om ) P_{|\fr}\otimes\eta$, 
for $P\in S(\fg^*)$, $\om\in\wedge (\fg^*)$, $\eta\in\Om (M)$; $\veps : 
\wedge(\fg^*)\to\bc$ being the standard augmentation map.
  \end{lemma}

  \begin{proof}  From the definition of $\Phi_M$, it suffices to prove 
that the following diagram is commutative:
  \beqn \tag{$*$} \begin{CD}  
\bigl( S(\fg^*)\otimes\Om (M_G)\bigr)^G @>{\bar{\psi}_M}>{\sim}> \bigl( 
W(\fg )\otimes\Om (M)\bigr)_R\\
@V{\hat{\cq}^{M_G}_G}VV  @VV{\cq^M_{G/R}}V \\
\bigl( U(\fg )\overset{\bullet}{\otimes}\Om (M_G)\bigr)^G 
@>{\psi_M}>{\sim}> \Ind_{G/R}(\Om (M)), 
  \end{CD}  \eeqn
where $\bar{\psi}_M$ is defined the same way as $\psi_M$. (Use the fact that
$\hat{\al}_M \circ \bar{\psi}_M=\al_M$.)   Considering the 
$G$-equivariant canonical projection (with $G$ acting on $G\times M$ via 
its left multiplication on the first factor) $G\times M\to M_G$, to prove 
the commutativity of ($*$), we can replace $M_G$ by the $G$-manifold 
$G\times M$.  From the definition of the various maps in ($*$), we can 
further assume that $M$ is the one point manifold $M^o$, i.e., we are 
reduced to prove the commutativity of ($*$) for $M_G$ replaced by $G$ 
(with $G$ acting on $G$ via the left multiplication).  Again using the 
definition of the various maps in ($*$), we are reduced to proving that 
  \beqn
\sig (x) = (-1)^n x, \; \text{ for all $(\wedge (\fg^*)\otimes\wedge 
(\fg^*))_{\text{hor}}$ of total degree $n$},
  \eeqn
where $\sig$ is the involution of $\wedge (\fg^*)\otimes\wedge (\fg^*)$ 
taking $\om\otimes\eta \mapsto (-1)^{\deg\om\cdot\deg\eta}\,\eta\otimes\om$, 
and `hor' is taken with respect to the standard tensor product action of 
$i_a$ ($a\in\fg$) on $\wedge (\fg^*)\otimes\wedge (\fg^*)$.

Clearly, $\del\om_1 := \om_1\otimes 1 - 1\otimes\om_1\in 
\bigl(\wedge (\fg^*)\otimes\wedge (\fg^*)\bigr)_{\text{hor}}$,
 for any $\om_1\in\fg^*$.  
Thus the subalgebra $\ca$ of $\wedge (\fg^*)\otimes\wedge (\fg^*)$ 
generated by $\{\del\om_1\}_{\om_1\in\fg^*}$ is contained in 
$\bigl(\wedge (\fg^*)\otimes\wedge (\fg^*)\bigr)_{\text{hor}}$.  Moreover, 
since the projection
  \[
\bigl(\wedge (\fg^*)\otimes\wedge (\fg^*)\bigr) \to \wedge (\fg^*), \; 
\om\otimes\eta \mapsto \veps (\eta )\om ,
  \]
induces an isomorphism
  \[
\bigl(\wedge (\fg^*)\otimes\wedge (\fg^*)\bigr)_{\text{hor}} \simeq \wedge 
(\fg^*),
  \]
we get that
  \beqn
\ca = \bigl(\wedge (\fg^*)\otimes\wedge (\fg^*)\bigr)_{\text{hor}} .
  \eeqn
Clearly (1) is satisfied for $\del\om_1$ and hence for each  
element of $\ca$.  Thus we get (1) by (2) and the lemma is proved.
  \end{proof}

\begin{remark} (1) An appropriate analogue of Theorems (2.2) and (2.5), Corollary (2.3) and Lemma (2.6) can be proved for any $R$-differential algebra $B$ replacing $\Om(M)$.

(2) Instead of defining $\Phi_M$ as in 2.4, we could 
have  (uniquely) defined $\Phi_M$ satisfying the above lemma.  But we find 
the original definition (as in 2.4) easier to work with.
\end{remark}

 \section{Cubic Dirac operator and results of\\ Huang-Pand\u{z}i\'c and 
Kostant}

We follow the notation and assumptions as in the beginning of Section 2.  
In particular, $(G,R)$ is a quadratic pair, i.e., $G$ is a  real Lie group with 
complexified Lie algebra $\fg$ and $R\subset G$ is a closed  
 subgroup with complexified Lie subalgebra $\fr\subset\fg$.  We 
assume that $\fg$ has a $G$-invariant nondegenerate symmetric bilinear 
form ${B_\fg} = \ip<\, ,\,>$ such that ${B_\fg}_{|\fr}$ is nondegenerate.

  We now  identify the differential of $\Ind_{G/R}(M)$ for $M$ a 
one point manifold with Kostant's cubic Dirac operator.

  \begin{lemma}  
Let $M$ be the one point manifold $M^o$.  Then $\Ind_{G/R}(M^o)$ can canonically be 
identified with the super-algebra $(U(\fg )\otimes\Cl (\fp ))^R$ (under the 
tensor product algebra structure).

Moreover, the differential $d$ on $\Ind_{G/R}(M)$ under the above 
identification is given by
  \[
d(x) = \ad \cd^\fp (x), 
  \]
where
  \[
\cd^\fp := \sum\, u_{p_\ell}\otimes c_{q_\ell} - 1\otimes \gamma_\fp ,
  \]
 $\{ p_\ell\}_\ell$ is any basis of $\fp$ and $\{ q_\ell\}_\ell$ is the 
dual basis with respect to the nondegenerate form ${B_\fg}_{ |\fp}$ and 
$\gamma_{\fp} \in \wed^3(\fp )\simeq \wed^3(\fp^*)$ is the Cartan form
  \[
\gamma_{\fp} (x,y,z) = \ip<x, [y,z]> , \;\text{ for }x,y,z\in \fp .
  \]
It is easy to see that $\cd^\fp$ is $R$-invariant, i.e., $\cd^\fp \in 
(U(\fg )\otimes\Cl (\fp ))^R$.
  \end{lemma}

  \begin{proof}  Let $\{ r_m\}$ be a basis of $\fr$ and let $\{ s_m\}$ be 
the dual basis of $\fr$ under $B_{\fg_{|\fr}}$.  Then, of course, $\{ 
r_m\}_m \cup \{ p_\ell\}_\ell$ is a basis of $\fg$ and $\{ s_m\}_m\cup\{ 
q_\ell\}_\ell$ is the dual basis of $\fg$.  Thus the element $\cd\in\cw 
(\fg )$ as in Section (1) is given by
  \[
\cd = \sum_m u_{r_m}\otimes c_{s_m} + \sum_\ell u_{p_\ell}\otimes 
c_{q_\ell} - 1\otimes \gamma_{\fg} .
  \]
Now
  \begin{align*}
\cw (\fg )_R &= (U(\fg )\otimes \Cl (\fg ))^R_{\hor R} \\
&\cong (U(\fg )\otimes \Cl (\fp ))^R . 
  \end{align*}
The differential $d$ in $\cw (\fg )$ is given by $dx = \ad \cd (x)$.  
Moreover, $d$ keeps the subspace $\cw (\fg )_R$ stable.  From this it 
is easy to see that, for $x\in \cw (\fg )_R$, $dx = \ad \cd^{\fp} (x)$.  
This proves the lemma.
  \end{proof}

  \begin{definition}  As in [Ko$_1$, \S1.5], the adjoint representation of 
$R$ on $\fp$ gives rise to the Lie algebra homomorphism
  \begin{align*}
\al : \fr &\to \Cl (\fp )^{\text{even}} \quad\text{ satisfying} \\
[\al (x), y] &= [x,y], \quad\text{ for $x\in\fr$ and $y\in\fp$}, 
  \end{align*}
where the bracket on the left side is commutation  in $\Cl (\fp )$.  Then, 
$\al$ is an $R$-module map under the adjoint actions. In particular, for 
$x_1,x_2\in\fr$, 
  \beqn
\al [x_1,x_2] = x_1\cdot \al (x_2). 
  \eeqn

Thus, we get an algebra homomorphism
  \[
\xi : U(\fr ) \to U(\fg )\otimes \Cl (\fp ), 
  \]
so that $\xi (x) = x\otimes 1 + 1\otimes\al(x)$, for $x\in\fr$.  It is 
easy to see that $\xi$ is injective.  Moreover, the earlier given 
$R$-module structure on $U(\fg )\otimes \Cl (\fp )$ 
(obtained from the adjoint action) is compatible with 
$\xi$. In particular,   for $x\in\fr$ and $a\in U(\fg )\otimes\Cl(\fp )$, 
  \beqn
x\cdot a = \xi (x)a - a\xi (x).
  \eeqn

Let $Z(G )$ (resp. $Z(R )$) be the subalgebra of invariants $U(\fg )^G$ (resp. $U(\fr 
)^R$).  Then, $Z(G )\otimes 1$ and $\xi (Z(R))$ are subalgebras of 
$(U(\fg )\otimes\Cl (\fp ))^R$.  Further, for $d = \ad \cd^{\fp}$, 
  \beqn
d_{|Z(G)\otimes 1} \equiv 0, 
  \eeqn
and
  \beqn
d_{|\xi (U(\fr ))} \equiv 0, 
  \eeqn
since $\xi (U(\fr ))$ commutes with any element in $(U(\fg )\otimes\Cl 
(\fp ))^R$ by (2).  Thus, by Theorem 
(2.5) and Lemma (3.1), we get algebra homomorphisms 
  \[
Z(G ) \to H\bigl( (U(\fg )\otimes\Cl (\fp ))^R,\, \ad \cd^{\fp}\bigr) = 
H\bigl( \Ind_{G/R}(\bc )\bigr) \overset{[\Phi_{M^o}]}{\tosim} Z(R ), 
  \]
where, as earlier,  $M^o$ is the one point manifold and the first map is induced from 
the map $z\mapsto z\otimes 1$.  Let $\eta_{R}$ be the composite algebra 
homomorphism 
  \[
\eta_{R} : Z(G) \to Z(R ).
  \]

Define a $R$-differential algebra homomorphism $F =F^\fg_\fr : W(\fr ) 
\to W(\fg )$ by 
  \[
F(\lam\otimes 1) = \bar{\lam}\otimes 1 - 1\otimes\delta (\lam ) \text{  and  
} F(1\otimes\lam ) = 1\otimes\bar{\lam}, \;\;\text{for }\lam\in \fr^*,
  \]
where $\bar{\lam}\in\fg^*$ is defined by $\lam_{|\fr} = \lam$ and 
$\bar{\lam}_{|\fp} \equiv 0$, and $\delta :\fr^* \to \wedge^2(\fg^*)$ 
is defined by
  \begin{alignat*}{2}
\delta (\lam )(y,z) &= \bar{\lam}([y,z]) \quad &&\text{ for $y,z\in\fp$}, 
\\
 &= 0  && \text{ if at least one of $y,z\in\fr$}.
  \end{alignat*}

Similarly, define a $R$-differential algebra homomorphism
 $
\cf = \cf^\fg_\fr : \cw (\fr ) \to \cw (\fg )$  by
$$\cf (x\otimes 1) = x\otimes 1 + 1\otimes\al (x), \;\text{ and}\,\,
\cf (1\otimes x) = 1\otimes x, \qquad\text{ for $x\in\fr$}.$$
Clearly,
  \beqn
  \cf_{|U(\fr )} = \xi .
  \eeqn

Then, interestingly, as proved by Alkseev-Meinrenken (private 
communication), we have:
  \beqn
\cq_\fg \circ F = \cf\circ\cq_\fr ,
  \eeqn
i.e., the following diagram is commutative:
  \[  \begin{CD}
W(\fr ) @>F>> W(\fg )\\
@VV{\cq_\fr}V @VV{\cq_\fg}V \\
\cw (\fr ) @>>{\cf}> \cw (\fg ).
  \end{CD}  \]

  \end{definition}

As a corollary of Theorem (2.5) we get the following.  This was 
conjectured by Vogan (actually Vogan conjectured a slightly weaker version) and proved by Huang-Pand\u{z}i\'c [HP, Theorem 3.4] 
in the case $R$ is a maximal compact subgroup of a connected reductive $G$.  The  case when $G$ and $R$ are connected and reductive was  
proved by Kostant [Ko$_3$, \S4.1].

  \begin{theorem}
For the differential $d := \ad \cd^{\fp}$ on $(U(\fg )\otimes\Cl (\fp ))^R$, 
  \beqn
\Ker d = \xi (Z(R )) \oplus \Imo d.
  \eeqn
In particular, $\xi (Z(R )) \simeq H\Bigl( (U(\fg )\otimes\Cl (\fp 
))^R, \ad \cd^{\fp}\Bigr)$.

  \end{theorem}

  \begin{proof}  We first prove that the composite map 
$\Phi_{M^o}\circ\xi$:
  \[
Z(R) \overset{\xi}{\longrightarrow} (U(\fg )\otimes\Cl (\fp ))^R
 \overset{\Phi_{M^o}}{\longrightarrow} Z(R )
  \]
is an isomorphism. (In fact, we will see during the proof of the next theorem that $\Phi_{M^o}\circ \xi$ is the identity map.) As earlier, let $\{U(\fr)^p\}_{p\geq 0}$ be the standard
filtration of the enveloping algebra $U(\fr)$ and let $Z(R)^p:=
U(\fr)^p \cap Z(R)$. By the definition of the map $\xi$, for $a\in Z(R)^p
\setminus Z(R)^{p-1}$, 
\[\xi(a)= a\otimes 1+ x,\]
for some $x\in \bigl(U(\fg)^{p-1} \otimes \Cl (\fp)\bigr)^R$. Thus, from the definition of the map $\Phi_{M^o}$ and the description of the isomorphism 
$\hat{Q}^{M_G}_G$ as in [AM, Proposition 6.5],
\[\Phi_{M^o}\circ\xi (a)= a \,\text{mod}\, Z(R)^{p-1}.\]
From this we see that $\Phi_{M^o}\circ\xi$ is an isomorphism. 

Since 
$\Phi_{M^o}$ induces an isomorphism in cohomology by Theorem (2.5), we get 
that the induced cohomology map 
  \[
[\xi ] : Z(R) \to H\bigl((U(\fg )\otimes\Cl (\fp ))^R, \ad 
\cd^{\fp}\bigr)
  \]
is an isomorphism.  From this of course (1) follows immediately.
\end{proof}

  \begin{theorem}  
The algebra homomorphism $\eta_R: Z(G )\to Z(R )$ 
is the unique homomorphism making the following diagram commutative:
  \beqn  \begin{CD}
Z(G ) @>{\eta_R}>> Z(R )\\
@V{H_{G}}VV @VV{H_{R}}V \\
S(\fg )^G @>{\beta_{R }}>> S(\fr )^R,  
  \end{CD} \tag{D} 
\eeqn
where $\beta_{R}$ is the restriction map under the identification $S(\fg 
) \simeq S(\fg^*)$, $S(\fr )\simeq S(\fr^*)$ induced by the bilinear form 
$B_{\fg}$, and $H_{G}$ (resp. $H_{R}$) is the inverse of the Duflo 
isomorphism of $\fg$ (resp. $\fr$) restricted to $S(\fg )^G$ (resp. 
 $ S(\fr )^R$). (Recall that for reductive $G$, $H_G$ coincides with the Harish-Chandra  isomorphism.) 

Thus, for $z\in Z(G )$, 
  \beqn
z\otimes 1 - \xi\bigl(\eta_R(z)\bigr) = \cd^{\fp}a_z + a_z\cd^{\fp}, 
  \eeqn
for some $a_z\in (U(\fg )\otimes\Cl(\fp )^{\text{\rm odd}})^R$.
  \end{theorem}

  \begin{proof}   With the notation as 
in the first diagram of Definition 2.4, for any $z\in Z(G )$, $z\otimes 
1\in (U (\fg )\obulltimes\Om (M^o_G))^G$ and, moreover, 
$$\psi_{M^o} (z\otimes 1) = z\otimes 1.$$
  Thus, by [AM, Proposition 6.5], 
  \[
\Phi_{M^o}(z\otimes 1) = D_{\fr}\bigl( (D^{-1}_{\fg}(z))_{|\fr}\bigr) ,
  \]
where $D_{\fg} : S(\fg )\to U(\fg )$ is the Duflo isomorphism under the 
identification $S(\fg )\simeq S(\fg^*)$, and similarly for $D_{\fr}$.
This gives that
  \[
\eta_R(z) = [\Phi_{M^o}](z\otimes 1) = D_{\fr}\circ \beta_{R}\circ 
(D_{\fg}^{-1})(z).
  \]

  From this 
the first part of the theorem follows.

We next prove that
  \beqn
\Phi_{M^o} \circ \xi_{|Z(R)} = I, 
  \eeqn
where $\xi : U(\fr ) \to U(\fg )\otimes\text{Cl}(\fp )$ is defined in 
\S3.2.  By Lemma (2.6), and the identities (3.2.5), (3.2.6), for $
x\in S(\fr^*)^R$, 
  \begin{align*}
\Phi_{M^o}\circ \xi\circ\cq_\fr (x) &= \hat{\cq}_R^{M^o} \circ 
\hat{\al}_{M^o}\circ F(x)\\
&= \cq_\fr \circ \hat{\al}_{M^o}\circ F(x), \;\text{ since } 
\cq_{\fr_{|S(\fr^*)^R} }=\hat{\cq}_R^{M^o} \\
&= \cq_\fr (x), \quad\; \text{ from the definition of $F$ and 
$\hat{\al}_{M^o}$}.
  \end{align*}
Since $\cq_{\fr_{|S(\fr^*)^R}} $ is an isomorphism onto $Z(R)$,  this 
proves (2).

From (2) we easily see that, for $z\in Z(G)$, 
$$\Phi_{M^o}(z\otimes 1) = 
\Phi_{M^o}\bigl(\xi (\eta_R(z))\bigr) = \eta_R(z),$$ and, moreover, 
by (3.2.3), (3.2.4), both of $z\otimes 1$ 
and $\xi 
(\eta_R(z))$ are cycles under $\text{ad } \cd^\fp$.  Thus they differ by a
coboundary, proving (1). This proves  
the theorem.

Alternatively, we can also obtain (1) in the special (but important) case 
where $G$ and $R$ are connected reductive groups (and $B_{\fg_{|\fr}}$ is 
nondegenerate) by using a result of Kostant as follows
.

By virtue of Theorem (3.3), define the map $\hat{\eta}_R: Z(G) \to Z(R)$ such that 
$z\otimes 1 - \xi (\hat{\eta}_{R}(z)) \in \Imo d$. Then it is easy to see 
that $\hat{\eta}_{R}$ is an algebra homomorphism. Moreover, by
 Kostant [Ko$_3$, Theorem 4.2] (generalizing the corresponding result 
 in the case when $R$ is a maximal 
compact subgroup of $G$ by Huang-Pand\u{z}i\'c [HP, Theorem 5.5]), 
$\hat{\eta}_{R}$ replacing ${\eta}_{R}$ also makes the diagram $(D)$ commutative. Thus $\hat{\eta}_{R}={\eta}_{R}$, proving (1).
  \end{proof}

  \begin{definition} Let $S$ be the space of spinors for 
  $\Cl(\fp )$, which is a  simple module of  $\Cl(\fp )$. Then, for any $U(\fg)$-module $V$, $V\otimes S$ is a  $U(\fg)\otimes \Cl(\fp )$-module under the componentwise action. In particular, the element $D^\fp \in 
(U(\fg)\otimes \Cl(\fp ))^\fr$ defined in Lemma (3.1) acts as a linear endomorphism 
$D^\fp_V$ on $V\otimes S$. 

Following Vogan, define the {\it Dirac cohomology}
$$ H_D(\fg, \fr;V)= \frac{\Ker\,D_V^\fp}{\Ker\,D_V^\fp \cap \Imo D_V^\fp}.$$
Since the element $D^\fp$ commutes with $\xi(U(\fr))$ (cf. $\S$3.2), both of $\Ker D_V^\fp$ and $\Imo D_V^\fp$ are $\fr$-submodules of  $V\otimes S$ via $\xi$. Thus 
$ H_D(\fg, \fr;V)$ has a canonical $\fr$-module structure.

Let $\chi : Z(\fg) \to \Bbb C$ be an algebra homomorphism, where  $ Z(\fg)$
is the center of $U(\fg)$. Recall that a $U(\fg)$-module $V$ is said to have {\it central character} $\chi$ if, for all $v\in V$ and $z\in Z(\fg)$,
$$ z\cdot v= \chi (z) v.$$ 

As an immediate consequence of Theorem (3.4), one gets the following corollary. Recall that this corollary was conjectured by Vogan in the case $R$ is a maximal compact subgroup of $G$ and proved in this case by Huang-Pand\u{z}i\'c 
[HP] and proved for general reductive pairs by Kostant [Ko$_3$].
\end{definition}

\begin{corollary}
Let $V$ be a $U(\fg)$-module with central character $\chi$. Then, for any 
$v\in  H_D(\fg, \fr;V)$ and $z\in Z(\fg)$,
$$\chi (z) v= \eta_R(z) v.$$

Of course,  the homomorphism $\eta_R: Z(\fg) \to Z(\fr)$ is completely determined from the diagram $(D)$ of Theorem (3.4).

Loosely speaking, the corollary asserts that the central character of any irreducible $\fr$-submodule of  $H_D(\fg, \fr;V)$ (if nonzero) determines the central character of $V$.
\end{corollary}
\begin{proof} We can clearly assume that $G$ and  $R$ are connected and thus
$Z(\fg)=Z(G)$ and $Z(\fr)=Z(R)$. By (3.4.1),
$$ z\otimes 1- \xi (\eta_R(z))= D^\fp a_z + a_z D^\fp,$$
for some $a_z \in (U(\fg)\otimes \Cl(\fp ))^R$. Thus,  for any $v_o\in \Ker D_V^\fp$, 
$$ (z\otimes 1) v_o- \eta_R(z) v_o\in \Imo D_V^\fp \cap \Ker  D_V^\fp,$$
since  $\eta_R(z) v_o\in  \Ker  D_V^\fp.$ Thus,
$\chi (z) v= \eta_R(z) v$ in $ H_D(\fg, \fr;V)$.
\end{proof}

Applying the definition of $\psi_M$ as in Theorem (2.2), for the case 
$R=G$ and a $G$-manifold $M$, interestingly we get an explicit expression 
for the inverse of the isomorphism $\Theta$. 

  \begin{lemma}
Take $R=G$ and a $G$-manifold $M$.  Then the inverse of the isomorphism
  \[\Theta=
\Theta_M : (\cw (\fg )\otimes\Om (M))_G \to 
(U(\fg )\obulltimes\Om (M))^G
  \]
(cf. Proposition 1.2) is given by the composition
  \[
(U(\fg )\obulltimes\Om (M))^G \overset{I\otimes\mu^*}{\tosim} (U(\fg 
)\obulltimes\Om (M_G))^G \overset{\psi_M}{\longrightarrow} (\cw (\fg 
)\otimes\Om (M))_G ,
  \]
where $\mu^* : \Om (M)\to\Om (M_G)$ is the $G$-module map induced from the 
$G$-equivariant smooth map
  \[
\mu : G \times^G M \to M, \qquad (g,m) \mapsto g\cdot m.
  \]
  \end{lemma}

  \begin{proof}  Since $\Theta_M$ is a vector space isomorphism, it suffices to prove that
  \[
\psi_M \circ (I\otimes\mu^*) \circ \Theta_M = I.
  \]
From the functoriality of $\Theta$, we have the following commutative diagram:
  \[  \begin{CD}
(\cw (\fg )\otimes\Om (M))_G @>{I\otimes\mu^*}>> (\cw (\fg )\otimes\Om 
(M_G))_G \\
 @VV{\Theta_M}V @VV{\Theta_{M_G}}V \\
(U(\fg )\obulltimes\Om (M))^G @>>{I\otimes\mu^*}> (U(\fg )\obulltimes\Om 
(M_G))^G .
  \end{CD}  \]

Take $a = \sum_i\, x_i\otimes\om_i \in (U(\fg )\obulltimes\Om (M))^G$.  
Then from the above commutative diagram:
  \begin{align*}
\Theta_M\circ\psi_M\circ \Theta_{M_G} \circ 
(I\otimes\mu^*)\circ\Theta^{-1}_M(a)  
&= \Theta_M \circ \psi_M \circ (I\otimes\mu^*)(a) \\
&= \Theta_M \circ \psi_M\Bigl( \sum_i\, x_i\otimes (\mu^*\om_i)\Bigr) \\
&= \sum_i\, x_i\otimes\bigl( (\mu^*\om_i)_{|1\times M}\bigr)\\
&= \sum_i\, x_i\otimes\om_i = a.
  \end{align*}
This gives
  \[
\psi_M \circ \Theta_{M_G} \circ (I\otimes\mu^*)\circ\Theta^{-1}_M = 
\Theta^{-1}_M .
  \]
Thus $\psi_M \circ \Theta_{M_G} \circ (I\otimes\mu^*) = I$ and 
hence, from the above commutative diagram again, $\psi_M\circ 
(I\otimes\mu^*)\circ\Theta_M = I$.  This proves the lemma.
  \end{proof}

\begin{remark} After an earlier version of this paper was distributed, E. Meinrenken informed me that he and Alekseev have obtained some results (unpublished) which  overlaps with our work. In particular, they also have obtained Theorems (3.3) and (3.4).

\end{remark}

   \end{document}